\documentclass[preprint]{elsarticle}  

\usepackage{graphicx,threeparttable,natbib}
\usepackage{amsmath,amsfonts,amssymb,booktabs}
\usepackage{graphicx,epsfig,color,url,lscape}
\usepackage{subcaption}
\usepackage{tikz}
\journal{}
\date{}
\begin{document}
\begin{frontmatter}
\title{Constraint Programming models for the parallel drone scheduling vehicle routing problem}

\author{Roberto Montemanni\corref{cor1}}
\ead{roberto.montemanni@unimore.it}

\author{Mauro Dell'Amico}
\ead{mauro.dellamico@unimore.it}

\cortext[cor1]{Corresponding author}

\address{Department of Sciences and Methods for Engineering, University of Modena and Reggio Emilia, Via Amendola 2, 42122 Reggio Emilia, Italy}

\begin{abstract}
Drones are currently seen as a viable way for improving the distribution of parcels in urban and rural environments, while working in coordination with traditional vehicles like trucks. In this paper we consider the parallel drone scheduling vehicle routing problem, where the service of a set of customers requiring a delivery is split between a fleet of trucks and a fleet of drones. We consider two variations of the problem. In the first one the problem is more theoretical, and the target is the minimization of the time required to complete the service and have all the vehicles back to the depot. In the second variant more realistic constraints involving operating costs, capacity limitation and workload balance, are considered, and the target is to minimize the total operational costs. We propose several constraint programming models to deal with the two problems. An experimental champaign on the instances previously adopted in the literature is presented to validate the new solving methods. The results show that on top of being a viable way to solve problems to optimality,  the models can also be used to  derive effective heuristic solutions and high-quality lower bounds for the optimal cost, if the execution is interrupted after its natural end. 
\end{abstract}
\end{frontmatter}

\section{Introduction}

In the last few years, drone technology has seen investments for billion of dollars, due to its potentials. Forbes \cite{for} defined such a phenomenon  the ``Drone Explosion''. Drones can be applied to many sectors, among which logistics, surveillance and disaster relief \cite{ottooptimization}. The most prominent application is probably logistics related to e-commerce, which has experienced an exponential growing in the last decades (Statista, \cite{sta}).  In \cite{BCG}, the authors 
forecast that autonomous vehicles will deliver about 80\% of all parcels in the following ten years. Several advantages associated with the use of aerial drones can be identified: they do not have to stick to the road network but can fly approximately straight line. Moreover, they are not affected by road traffic congestions. The technology might lead to innovative solutions of interest for the companies (operational costs reduction), for the customers (faster deliveries) and for the whole society (sustainability). In this work we will focus on operational delivery strategies with a mixed fleet using both trucks and drones.

The seminal work of Murray and Chu \cite{murray2015flying} pioneered a new routing problem in which a truck and a drone collaborate to make deliveries. From an operations research perspective, the authors present two new prototypical variants expanding from the traditional Traveling Salesman Problem (TSP) called the Flying Sidekick TSP (FSTSP) and the Parallel Drone Scheduling TSP (PDSTSP). In both cases a truck and drones collaborate to deliver parcels, the difference being however that in the former model drones can be launched and collected from the truck during its tour, while in the latter one drones are operated directly from the central depot, while the truck executes a traditional delivery tour. In the remainder of the paper we will focus on the latter problem, addressing the interested reader, for example, to \cite{amicobb} for full details and some solution strategies for the FSTSP.

More formally, in the PDSTSP there is a truck that can leave the depot, serve a set of customers, and return to the depot, and a set of drones, each of which in the meantime can leave the depot, serve a customer, and return to the depot before serving other customers. Not all the customer can be served by the drones, either due to their location or the characteristic of their parcel. The objective of the problem is to minimize the completion time of the last vehicle returning to the depot (or a cost function related to this), while serving all the customers. 

A first Mixed Integer Linear Programming (MILP) model for the PDSTSP and some simple heuristic methods are proposed in \cite{murray2015flying}. A more refined mixed integer programming model and the first metaheuristic method, based on a two steps strategy, embedding a dynamic programming-based component, are discussed in \cite{mbiadou2018iterative}. A similar two steps approach, but based on matheuristics concepts is presented in \cite{DMN}. A hybrid ant colony optimization metaheuristic is discussed in 
\cite{dinh2022}. In \cite{md23} a constraint programming approach is discussed, which is able to solve to optimality all the benchmark instances previously adopted in the literature for both exact and heuristic methods. An improved variable neighbour search metaheuristic is discussed in \cite{lei2022}. More recently, in \cite{HA} another exact approach based on branch-and-cut was proposed, together with some new benchmark problems.

Several PDSTSP variants are also introduced and studied in the literature. We refer the interested reader to \cite{ottooptimization} and \cite{pasha2022} for a complete survey. In the following we review only the extensions of the original problem relevant to the present study, where multiple trucks are employed out of a same depot, and more realistic constraints involving load balancing, capacity and costs are eventually considered. 

The recent work  \cite{mbiadou2022} discusses the \emph{Parallel Drone Scheduling Multiple Traveling Salesman Problem} what we will here refer to as the MT-PDSVRP (\emph{ Min-Time Parallel Drone Scheduling Vehicle Routing Problem}), which is a straightforward extension of the PDSTSP where multiple trucks are employed and the target is to minimize the time required to complete the delivery to the last customer serviced and go back to the depot. The authors propose a hybrid metaheuristic algorithm along the line of the method previously introduced in \cite{mbiadou2018iterative}, a mixed integer linear model and a branch-and-cut approach working on such a model. What is basically the same problem is also introduced at the same time in \cite{raj2021}, where the authors propose three mixed integer linear programming models, one of which is arc-based and the other two are set covering-based, together with a branch-and-price approach based on one of the set covering-based models. A heuristic version of the branch-and-cut method is also discussed, targeting larger instances. A more realistic variation of the PDSTSP, which we will refer to as the MC-PDSVRP (\emph{Min-Cost Parallel Drone Scheduling Vehicle Routing Problem}) is introduced in \cite{nguyen2022}. In this version of the problem concepts such as capacity, load balancing and decoupling of costs and times are taken into account. A formal definition of the problem will be provided in Section \ref{vrpc}. The authors propose a mixed integer linear programming model and a ruin\&recreate metaheuristic for the problem. 

{In this paper we aim to explore the potential of constraint programming on PDSVRPs, in the hope of exploiting the recent advances of solvers in dealing with TSP and VRP problems and we present tailored models for the different PDSVRP considered. Experimental results show that several best-known upper and lower bounds can be obtained by the new models.}

{The rest of the paper is organized as follows. Section \ref{vrpt} is devoted to the the MT-PDSTSP. The problem is formally defined and two Constraint Programming models are presented. Symmetrically, in Section \ref{vrpc} the MC-PDSVRP is discussed, again presenting Constraint Programming models. Experimental results for both the models are discussed in Section \ref{sec:experiments}, while Section \ref{conc} contains some conclusions and ideas for future work.}

\section{The Min-Time Parallel Drone Scheduling Vehicle Routing Problem}\label{vrpt}
\begin{sloppypar}
The MT-PDSVRP can be represented on a complete directed graph $G~=~(V, A)$, where the node set  $V = \{0, 1, ..., n\}$ represents the depot (node $0$) and a set of customers  $C = \{1, ..., n\}$ to be serviced.
A set $T$ of homogeneous trucks and a set $D$ of homogeneous drones are available to deliver parcels to the customers. Each truck starts from the depot $0$, visits a subset of the customers, and returns back to the depot, operating a single route. The drones operate back and forth trips from the depot to customers, delivering one parcel per trip and operating multiple routes if necessary. Not all the customers can be served by a drone, due to the weight of the parcel, an excessive distance of the customer location from the depot, or eventual terrain obstacles such as hills or areas with high-rise buildings. Let $C^D \subseteq C$ denotes the set of customers that can be served by drones. These customers are referred to as \emph{drone-eligible} in the remainder of the paper. The travel time incurred by a truck to go from node $i$ to node $j$ is denoted as $t_{ij}^T$, while the time required by a drone to
serve a customer $i$ (back and forth) is denoted as $t_{i}^U$. The trucks and the drones start from the depot at time 0, and the objective of the MT-PDSVRP is to minimize the time required to complete all the deliveries and to have all the trucks and all the drones back at the depot. Note that since truck and drones work in parallel, the objective function translates into minimizing the  time required by the vehicle of the fleet with the longest total operational time.
\end{sloppypar}

\begin{figure}[h]
	
	\begin{subfigure}{.45\textwidth}
		
		\begin{tikzpicture}[node distance={2cm}, main/.style = {draw, circle}]
			\node[main,minimum size=0.75cm,double] (0) {0}; 
			\node[main,minimum size=0.75cm] (1) [left of=0] {1};
			\node[main,minimum size=0.75cm] (2) [above left of=0] {2};
			\node[main,minimum size=0.75cm] (3) [above right of=2] {3};
			\node[main,minimum size=0.75cm] (4) [right of=3] {4};
			\node[main,minimum size=0.75cm] (6) [below right of=0] {6};
			\node[main,minimum size=0.75cm] (7) [below left of=6] {7};	
			\node[main,minimum size=0.75cm] (8) [left of=7] {8};
			\node[main,minimum size=0.75cm] (5) [right of=0] {5};		
		\end{tikzpicture}
	\end{subfigure}\hspace{20mm}%
\begin{subfigure}{.45\textwidth}
		
		\begin{tikzpicture}[node distance={2cm}, main/.style = {draw, circle}]
			\node[main,minimum size=0.75cm,double] (0) {0}; 
			\node[main,minimum size=0.75cm] (1) [left of=0] {1};
			\node[main,minimum size=0.75cm] (2) [above left of=0] {2};
			\node[main,minimum size=0.75cm] (3) [above right of=2] {3};
			\node[main,minimum size=0.75cm] (4) [right of=3] {4};
			\node[main,minimum size=0.75cm] (6) [below right of=0] {6};
			\node[main,minimum size=0.75cm] (7) [below left of=6] {7};	
			\node[main,minimum size=0.75cm] (8) [left of=7] {8};
			\node[main,minimum size=0.75cm] (5) [right of=0] {5};			
			\path [draw,->,>=stealth, line width=1.2] (0) ->  (2);
			\path [draw,->,>=stealth, line width=1.2] (2) -> (3);
			\path [draw,->,>=stealth, line width=1.2] (3) -> (0);
			\path [draw,->,>=stealth, line width=1.2, color=blue] (0) ->  (6);
			\path [draw,->,>=stealth, line width=1.2, color=blue] (6) -> (7);
			\path [draw,->,>=stealth, line width=1.2, color=blue] (7) -> (0);
			\path[->,>=stealth, line width=1.2, color=red,dashed] (0) edge [bend right=15] node[auto] {} (1);
			\path[->,>=stealth, line width=1.2, color=red,dashed] (1) edge [bend right=15] node[auto] {} (0);
			\path[->,>=stealth, line width=1.2, color=red,dashed] (0) edge [bend right=15] node {} (8);
			\path[->,>=stealth, line width=1.2, color=red,dashed] (8) edge [bend right=15] node {} (0);
			\path[->,>=stealth, line width=1.2, color=green,dashed] (0) edge [bend right=15] node {} (5);
			\path[->,>=stealth, line width=1.2, color=green,dashed] (5) edge [bend right=15] node {} (0);
			\path[->,>=stealth, line width=1.2, color=green,dashed] (0) edge [bend right=15] node {} (4);
			\path[->,>=stealth, line width=1.2, color=green,dashed] (4) edge [bend right=15] node {} (0);
		\end{tikzpicture}
	\end{subfigure}\hspace{9mm}%
	\caption{{{An example} 
 of an instance is provided on the left; we assume that two trucks and two drones are available (travel times are omitted for the sake of simplicity). A solution of the MT-PDSVRP  is provided on the right. The black and blue arcs represent the tour of the two truck. The first visits nodes 2 and 3 before going back to the depot, the second visits nodes 6 and 7 and goes back to the depot. The dashed red arcs depict the missions of the first drone (that serves nodes 1 and 8), while the dashed green arcs depict the missions of the second drone (that serves nodes 4 and 5)}.}
	\label{fig-pcmca-example}
\end{figure}
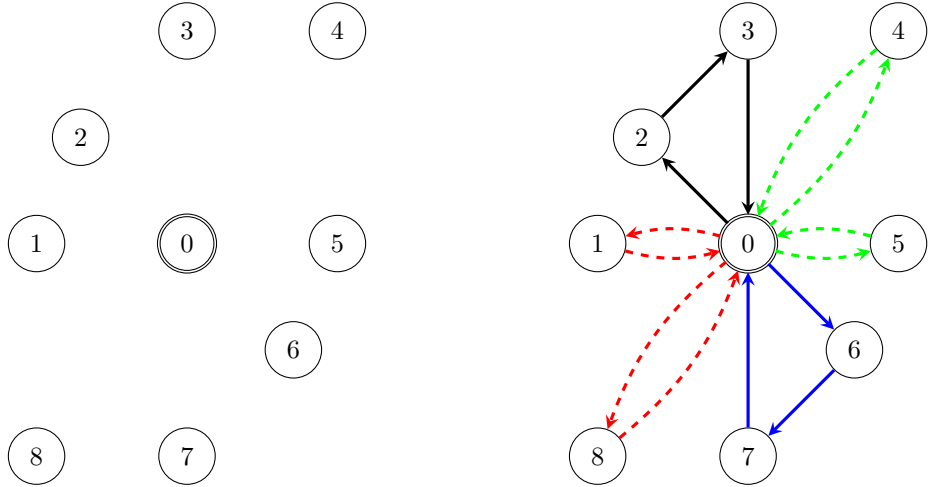

An example of a solution of the MT-PDSVRP for a small instance is provided in {Figure}
 \ref{fig-pcmca-example}.

\subsection{Constraint programming models} \label{cpt}

We propose two alternative models that exploit different functions made available by modern Constraint Programming solvers to model the problem. In the first model the truck tours are represented as separated non-overlapping entities, while in the second they are modelled similarly to the classic giant tour representation \cite{old}.

\subsubsection{MT-3IDX: A model based on separated truck tours}\label{t3}
The variables of the model that will be referred to as MT-3IDX are as follows. The binary variable $z_{kij}$, with $k \in T$ and $i,j \in V$ takes value 1  (or equivalently \emph{True}) if node $i$ is visited right before node $j$ in the tour of truck $k$, and value $0$  (or equivalently \emph{False}) otherwise. The binary variable $x_{di}$, with $d \in D$ and $i \in C^D$, takes value $1$ is the customer $i$ is visited by the drone $d$, and value $0$ otherwise. Finally, the variable $\alpha$ is a continuous variable introduced to implement the $\min\max$ objective function. The model is as follows.
\begin{align} 
  MT-&3IDX:  \ \ \ \min   \alpha & \label{1a}\\
 s.t. \ \ \ &	\text{Circuit($z_{kij}; i, j \in V; i \neq 0 \lor j \neq 0$)} & k \in T \label{2a}\\
      & \sum_{d\in D} x_{di} +  \sum_{k\in T} \neg z_{kii} = 1    & i \in C^D \label{3a}\\
      & \sum_{k \in T} \sum_{i \ in V} \neg z_{kij}  = 1   \hspace{-2cm} & j \in C \setminus C^D  \label{33a}\\
       & \alpha \geq \sum_{i \in V} \sum_{j \in V, i \ne j} t^T_{ij} z_{kij} & k \in T   \label{4a}\\
      & \alpha \geq \sum_{i \in C^D} t^D_{i} x_{di} &d \in D \label{5a}\\
      & z_{kij} \in \{0;1\} & k \in T; i, j \in V \label{6a} \\
	& x_{di} \in \{0; 1\} & d\in D, i \in C^D \label{7a}
\end{align}
The objective function (\ref{1a}) minimizes $\alpha$, that will be assigned the time taken by the latest vehicle to complete its tasks.
{Constrains (\ref{2a}) assures that the $z$ variables associated to each truck $k$ take values in such a way to form a valid tour, eventually with self-loops for the variables associated to customers visited by the drones. This is imposed through the use of the ``\emph{Circuit}'' statement.} In the logic of such a command, if a customer $i$ is not visited by the a truck $k$ then $z_{kii}$ is set to 1.
Constrains (\ref{3a}) imposes that if a customer $i \in C^D$ is not visited by any truck, then one of the drones in $D$ must visit it. Note that the negation operator ``$\neg$'' is used, and that if $z_{kii}=1 \ \forall k \in T$, then customer $i$ is not visited by any truck.
Constraints (\ref{33a}) imposes that if a customer $i \in C \setminus C^D$ then one truck in $T$ must visit it. 
Inequality (\ref{4a}) constraints $\alpha$ to be at least as large as the length of the tour of each truck.
Constrains (\ref{5a}) set $\alpha$ to be at least as large as the time spent by each drone $d$ to execute the tasks assigned to it.
Finally, constrains (\ref{6a}) and (\ref{7a}) define the variables domain.  

This model  presents the drawback of having many variables, and these variables also have to be constrained to take mutually-consistent values. This might affect the quality of the lower bounds provided by the relaxations of the model.

\subsubsection{MT-2IDX: A model based on a giant truck tour}\label{t4}
In this section we present a model with less variables than that discussed in Section \ref{tc}. This is achieved by exploiting the ``\emph{MultipleCircuit}'' statement present in modern Constraint Programming tools \cite{ortools}. 

The meaning of variables $x$ remain the same as in Section \ref{tc}, while $z$ variables are substituted by a new set of variable $y_{ij}$, with $i,j \in V$, {that when $i \neq j$ takes value 1 if node $i$ is visited right before node $j$ in one of the truck tours, and value $0$ otherwise.  When $i=j$, $y_{ii}$ takes value 1 if node $i$ is visited by a drone, 0 otherwise.} A new set of continuous variables $\gamma_i$, $i \in V$ is introduced to count the operating time of each truck $k \in T$.  The resulting  model is as follows.
\begin{align} 
 MT-&2IDX: \ \ \ \min   \alpha & \label{1bab}\\
 s.t. \ \ \ &	\text{MultipleCircuit($y_{ij}, i,j \in V;  i \neq 0 \lor j \neq 0$)} \hspace{-3cm}& \label{2bab}\\
	& \sum_{j \in V \setminus \{0\}} y_{0j} \leq |T|  \hspace{-3cm}& \label{3bab}\\
      & \sum_{d\in D} x_{di} +  \neg y_{ii}  = 1   \hspace{-3cm}  & i \in C \label{4bab}\\
      & \gamma_0=0 & \label{9bab}\\
      & y_{ij} \implies   \gamma_j = \gamma_i + t^T_{ij}  & i \in V;  j \in V \setminus \{0\} \label{10bab}\\
      & y_{i0} \implies    \alpha \geq \gamma_i + t^T_{i0} & i \in V \setminus \{0\}  \label{11bab}\\
       & \alpha \geq \sum_{j \in C^D} t^D_j x_{dj}   \hspace{-3cm}& d \in D  \label{5bab}\\
      & y_{ij} \in \{0;1\} & i, j \in V \label{12bab}\\
	& x_{di} \in \{0; 1\} & d\in D, i \in C^D \label{13bab}\\
      & \gamma_i \ge 0 & i \in V \label{14bab}
\end{align}
The objective function (\ref{1bab}) minimizes $\alpha$, that will be assigned the time taken by the latest vehicle to complete its tasks.
Constrains (\ref{2bab}) assures that the $y$ variables take values in such a way to form a set of valid truck tours, and eventually self-loops for the variables associated to customers visited by the drones. This is imposed through the use of the ``\emph{MultipleCircuit}'' statement.
Constrain (\ref{4bab}) states that the number of truck tours (outgoing arcs from the depot $0$) is limited to $|T|$.
Constrains (\ref{3bab}) imposes that if each customer $i \in C$ must be visited either by a truck ($y_{ii}=0$) or by one of the drones in $D$.
Constraint (\ref{9bab}) initializes to 0 the time counter at the depot.
Constraints (\ref{10bab}) increase the operational time of each truck (by setting the time each customer $i$ is visited). 
{Inequalities (\ref{11bab}) force $\alpha$ to be at least as large as the last return time of a truck to the depot.} 
Constrains (\ref{5bab}) counts the total operating time of each drone $d$ and imposes that $\alpha$ must be greater than it.
Finally, constrains (\ref{12bab}), (\ref{13bab}) and (\ref{14bab}) define the domain of the variables.  

{With respect to the model discussed in \ref{tc},  MT-2IDX has the advantage of having a few less variables, but on the other hand it requires more complex constraints to account for the length of each single truck tour.}
  
\section{The Min-Cost Parallel Drone Scheduling Vehicle Routing Problem}\label{vrpc}
The  MC-PDSVRP is a modification of the  MT-PDSVRP described in Section \ref{vrpt}, where more realistic components are considered, {leading to a model with more complex objective function and more constraints.}. In particular, the following elements are added to concepts already analyzed for the  MT-PDSVRP and have to be taken into account:
{A transportation cost $c_{ij}^T$ is incurred by a truck to travel on the edge $(i,j) \in E$, while $c^D_i$ is the transportation costs incurred by each drone to complete a mission to customer $i \in C^D$. These costs are used to define a new objective function to be minimized. The delivery request from each customer $i \in C$ is associated with a parcel weight $w_i$. Each truck has to respect some capacity constraints, being limited to transport at most $Q^T$ units of weight in its route. Finally, there are upper bounds $\tau^T$ and $\tau^D$ on the maximum working time for each truck and drone, respectively. Note that these latter constraints are normally used to introduce load-balancing concepts into the optimization. }

\subsection{Constraint programming models} \label{cpc}
Analogously to what seen in Section \ref{cpt} for the  MT-PDSVRP, two models will be presented. 

\subsubsection{MC-3IDX: A model based on separated truck tours}\label{tc}

{The variables used in this model are the same used for the model MT-3IDX in Section \ref{t3}.}

\begin{align} 
 MC-&3IDX:   \ \ \ \min  \sum_{k \in T}\sum_{i \in V}\sum_{j \in V, j \ne i} c^T_{ij}z_{kij} + \sum_{k \in D}\sum_{i \in V} c^D_{i} x_{ki}\hspace{-4cm}& \label{1d}\\
  s.t. \ \ \ &	\text{Circuit($z_{kij}; i,j \in V$)} \hspace{-2cm}& k \in T\label{2d}\\
      &  \sum_{d \in D} x_{dj} + \sum_{k \in T} \sum_{i \ in V} \neg z_{kij} = 1   \hspace{-2cm} & j \in C^D  \label{3d}\\
      & \sum_{k \in T} \sum_{i \ in V} \neg z_{kij}  = 1   \hspace{-2cm} & j \in C \setminus C^D  \label{5d}\\
       & \sum_{i \in V} \sum_{j \in V, i \ne j} t^T_{ij} z_{kij} \leq \tau^T & k \in T  \label{7d}\\
       & \sum_{j \in C^D} t^D_j x_{dj} \leq \tau^D & d \in D  \label{6d}\\
       & \sum_{j \in C} w_{j} z_{kij} \leq Q^T & k \in T  \label{8d}\\
      & z_{kij} \in \{0;1\} & k \in T, i, j \in V \label{13d}\\
     	& x_{di} \in \{0; 1\} & d\in D, i \in C^D \label{14d}
\end{align}

{The objective function (\ref{1d}) minimizes the total cost, given by sum of the costs of each truck and drone, incurred to service all the customers.}
Constrains (\ref{2d}) assures that the $z$ variables associated to each truck $k$ take values in such a way to form a valid tour, eventually self-loops for the variables associated to customers visited by the drones.
Constraints (\ref{3d}) imposes that if a customer $i \in C^D$ is not visited by any truck, then one of the drones in $D$ must visit it. 
Constraints (\ref{5d}) imposes that if a customer $i \in C \setminus C^D$ then one truck in $T$ must visit it. 
Inequalities (\ref{7d}) make sure that each truck does not exceed its maximum working time. 
Inequalities (\ref{6d}) are the constraints on maximum working time for the drones. 
Constraints (\ref{8d}) model the capacity constraints for the trucks.
Finally, constraints (\ref{13d}) and (\ref{14d}) define the variables domain.  

\subsubsection{MC-2IDX: A model based on a giant truck tour}

{The variables used in this model are the same used for the model MT-2IDX in Section \ref{t4}, with the addition of the new set of continuous variables $\beta_i$, $i \in V$ which is introduced to count the weight carried by each truck $k \in T$.}

\begin{align} 
 MC-&2IDX:   \ \ \ \min \sum_{i \in V}\sum_{j \in V, j \ne i} c^T_{ij}y_{ij} + \sum_{d \in D}\sum_{i \in V} c^D_{i} x_{di}\hspace{-4cm}& \label{1ba}\\
 s.t. \ \ \ &	\text{MultipleCircuit($y_{ij}, i,j \in V;  i \neq 0 \lor j \neq 0$)} \hspace{-3cm}& \label{2ba}\\
      & \sum_{d\in D} x_{di} +  \neg y_{ii}  = 1   \hspace{-3cm}  & i \in C \label{3ba}\\
	& \sum_{j \in V \setminus \{0\}} y_{0j} \leq |T|  \hspace{-3cm}& \label{4ba}\\
       & \sum_{j \in C^D} t^D_j x_{dj} \leq \tau^D  \hspace{-3cm}& d \in D  \label{5ba}\\
      & \beta_0=0 & \label{6ba}\\
      &y_{ij} \implies  \beta_j = \beta_i + w_j   \!\!& i \in V;  j \in V\setminus \{0\}  \label{7ba}\\
      &  \beta_i \leq Q^T  & i \in V \setminus \{0\} \label{8ba}\\
      & \gamma_0=0 & \label{9ba}\\
      & y_{ij} \implies   \gamma_j = \gamma_i + t^T_{ij}  & i \in V;  j \in V \setminus \{0\} \label{10ba}\\
      & y_{i0} \implies   \gamma_i + t^T_{i0} \leq \tau^T & i \in V \setminus \{0\}  \label{11ba}\\
      & y_{ij} \in \{0;1\} & i, j \in V \label{12ba}\\
	& x_{di} \in \{0; 1\} & d\in D, i \in C^D \label{13ba}\\
      & \beta_i, \gamma_i \ge 0 & i \in V \label{14ba}
\end{align}
The objective function (\ref{1ba}) minimizes the total cost, given by sum of the costs of each truck and drone, incurred to service all the customers.
Constrains (\ref{2ba}) assures that the $y$ variables take values in such a way to form a set of valid truck tours, and eventually self-loops for the variables associated to customers visited by the drones. 
Constrains (\ref{3ba}) imposes that if each customer $i \in C$ must be visited either by a truck ($y_{ii}=0$) or by one of the drones in $D$.
Constrain (\ref{4ba}) states that the number of truck tours (outgoing arcs from the depot $0$) is limited to $|T|$.
Constrains (\ref{5ba}) limits the total operating time of each drone $d$ to the maximum allowed value $\tau^D$.
Constraint (\ref{6ba}) initializes to 0 the weight counter at the depot.
Constraints (\ref{7ba}) increase the weight transported by each truck. 
Inequalities (\ref{8ba}) impose that the incremental weight $\beta_i$ at each customer $i \in V$ can never exceed the maximum capacity $Q^T$ of each truck. 
Constraint (\ref{9ba}) initializes to 0 the time counter at the depot.
Constraints (\ref{10ba}) increase the operational time of each truck (by setting the time each customer $i$ is visited). 
Inequalities (\ref{11ba}) impose that the incremental time $\gamma_i$ at each customer $i \in V$, plus the time required to the truck to go back to the depot, can never exceed the maximum operating time $\tau^T$ of each truck. 
Finally, constrains (\ref{12ba}), (\ref{13ba}) and (\ref{14ba}) define the variables domain.  
 
\section{Computational experiments} \label{sec:experiments}
The constraint programming models described in Sections \ref{vrpt} and \ref{vrpc} have been implemented in Python 3.9 and solved via the CP-SAT solver of Google OR-Tools 9.5.2237 \cite{ortools}. The experiments have been run on a laptop computer equipped with 32 GB of RAM and an Intel Core i7 12700F CPU with 12 cores (8 with a maximum frequency of 4.9 GHz and 4 with a maximum frequency of 3.6 GHz).

The outcome of the experimental champaign we run is discussed in the remainder of this section, and is organized according to the different problems attacked. All the tables present the information of the instances and for each method considered the given maximum computation time (in the label of the columns), although the hardware is different in each paper. For each instance/method combination,  the cost of the best heuristic solution retrieved, and eventually the lower bound found. A dash means no result was retrieved in the given time {(or, for some of the methods we compare with, the experiment was not attempted)}. {Every time one of the new models we propose improve a best-known bound, the relative entry is in bold. Analogously, every time one of the new models we propose does not match or improve a best-known bound, the relative entry is in italic. Finally, proven exact solutions, retrieved by any method, are marked with an asterisk in the tables.} All the new best-known solutions retrieved are available upon request to the authors.

\subsection{MT-PDSVRP} \label{rvrpt}	
{The results are subdivided based on the source of the instances in the following subsections. }

\subsubsection{Instances from Mbiadou Saleu et al. \cite{mbiadou2022}} 

A first set of benchmarks for the MT-PDSVRP was created in \cite{mbiadou2022} starting from classic instances for the \emph{Capacitated Vehicle Routing Problem}. A total of 20 instances with a number of customers ranging between 50 and 199 has been obtained. We refer the interested reader to \cite{mbiadou2022}  for full details about the elements of the new instances and the complete sources for the instances. The results are summarized in Table \ref{tf}, where the results are those obtained by a branch and cut (BC) approach running on a MILP model and of the best of nine variations of a hybrid metaheuristc approach (HM). The BC solver is run on a computer with an Intel Xeon CPU E5-2670 CPU with 2x8 cores running at 2.6 GHz and 62.5 GB of RAM, with a maximum computation time of 10800 seconds, while the HM variations have been run on a computer with an Intel core i5-6200U CPU running at 2.4 GHz and 8 GB of RAM, with a maximum computation time of 1000 for each variation.

\begin{table}[h!]														
\caption{MT-PDSVRP - Instances originated from CVRPLIB in Mbiadou Saleu et al. \cite{mbiadou2022}.}													
\label{tf}																	
{		
\hspace{-2cm}							
\footnotesize
 \begin{tabular}{| l c c | cc | c | cc|cc|}		\hline															
\multicolumn{3}{|c|}{Instances}	&		\multicolumn{3}{|c|}{Mbiadou Saleu et al.  \cite{mbiadou2022}}		&		\multicolumn{2}{|c|}{MT-3IDX}&		\multicolumn{2}{|c|}{MT-2IDX}	\\
\multicolumn{3}{|c|}{}	&		\multicolumn{2}{|c|}{ BC (10800 sec)}	&		HM	(9x1000 sec)&		\multicolumn{2}{|c|}{(10800 sec)}&		\multicolumn{2}{|c|}{(10800 sec)}	\\
Name		&		$|T|$		&		$|D|$		&		LB	&		UB		&	UB	&	LB	&	UB&	LB	&	UB	 	\\ \hline
CMT1	&	3	&	2	&	145.87	&	188.00	&	166.00	&	\textbf{	158.27	}	&	\emph{	168.00	}	&	\emph{	103.25	}	&	\emph{	174.00	}	\\
CMT2	&	5	&	5	&	101.66	&	3630.86	&	130.23	&	\textbf{	108.00	}	&	\emph{	158.00	}	&	\emph{	65.33	}	&	\emph{	144.00	}	\\
CMT3	&	4	&	4	&	161.07	&	4537.11	&	184.00	&	\textbf{	165.98	}	&	\emph{	222.00	}	&	\emph{	98.07	}	&	\emph{	210.00	}	\\
CMT4	&	6	&	6	&	-	&	-	&	160.38	&	\textbf{	125.39	}	&	\emph{	288.00	}	&	\emph{	72.73	}	&	\emph{	180.00	}	\\
CMT5	&	9	&	8	&	-	&	-	&	138.00	&	\textbf{	82.84	}	&	\emph{	312.00	}	&	\emph{	55.67	}	&	\emph{	158.00	}	\\
E-n51-k5	&	3	&	2	&	145.87	&	188.00	&	168.00	&	\textbf{	159.34	}	&	\textbf{	166.00	}	&	\emph{	103.22	}	&	\emph{	174.00	}	\\
E-n76-k8	&	4	&	4	&	127.05	&	2975.51	&	154.00	&	\textbf{	132.49	}	&	\emph{	180.00	}	&	\emph{	81.13	}	&	\emph{	164.00	}	\\
E-n101-k8	&	4	&	4	&	161.07	&	4537.11	&	184.00	&	\textbf{	165.64	}	&	\emph{	244.00	}	&	\emph{	98.03	}	&	\emph{	224.00	}	\\
M-n151-k12	&	6	&	6	&	-	&	-	&	154.00	&	\textbf{	126.57	}	&	\emph{	270.00	}	&	\emph{	72.78	}	&	\emph{	188.00	}	\\
M-n200-k16	&	8	&	8	&	-	&	-	&	144.00	&	\textbf{	95.21	}	&	\emph{	302.00	}	&	\emph{	59.19	}	&	\emph{	166.00	}	\\
P-n51-k10	&	5	&	5	&	81.35	&	230.00	&	111.07	&	\textbf{	90.13	}	&	\emph{	118.00	}	&	\emph{	52.60	}	&	\emph{	124.00	}	\\
P-n55-k7	&	4	&	3	&	101.49	&	308.00	&	126.00	&	\textbf{	111.00	}	&	\emph{	136.00	}	&	\emph{	70.40	}	&	\emph{	130.00	}	\\
P-n60-k10	&	5	&	5	&	84.30	&	246.00	&	114.00	&	\textbf{	95.00	}	&	\emph{	124.00	}	&	\emph{	54.89	}	&	\emph{	120.00	}	\\
P-n65-k10	&	5	&	5	&	94.66	&	580.00	&	126.00	&	\textbf{	105.41	}	&	\emph{	144.00	}	&	\emph{	61.46	}	&	\emph{	138.00	}	\\
P-n70-k10	&	5	&	5	&	99.42	&	3166.25	&	128.00	&	\textbf{	108.00	}	&	\emph{	164.00	}	&	\emph{	64.34	}	&	\emph{	150.00	}	\\
P-n76-k5	&	3	&	2	&	181.36	&	280.00	&	200.00	&	\textbf{	189.91	}	&	\emph{	206.00	}	&	\emph{	129.06	}	&	\emph{	214.00	}	\\
P-n101-k4	&	2	&	2	&	321.80	&	4725.47	&	342.00	&	\textbf{	329.90	}	&	\textbf{	340.00	}	&	\emph{	198.00	}	&	\emph{	362.00	}	\\
X-n110-k13	&	7	&	6	&	-	&	-	&	1864.00	&	\textbf{	1365.68	}	&	\emph{	2702.00	}	&	\emph{	766.70	}	&	\emph{	2048.00	}	\\
X-n115-k10	&	5	&	5	&	-	&	-	&	2258.00	&	\textbf{	1850.38	}	&	\emph{	2920.00	}	&	\emph{	1031.21	}	&	\emph{	2488.00	}	\\
X-n139-k10	&	5	&	5	&	-	&	-	&	2492.00	&	\textbf{	 -	}	&	\emph{	 -	}	&	\textbf{	1118.67	}	&	\emph{	2744.00	}	\\
\hline
\end{tabular}																	
}																	
\end{table}		

The results presented in Table \ref{tf} indicate that the new CP-based approaches are competitive with state-of-the-art results. In particular, new improved lower bounds were provided for all the instances considered, and 2 new best-known heuristic solutions were also retrieved. The comparison between the CP models indicate that the model MT-3IDX discussed in Section \ref{t3}, and characterized by the set of variables $z$ with 3 indices, performs better. In particular, it provides substantially tighter lower bounds. However, this model is not capable of handling {the last instance, which by its nature allows a large number of possible drone missions, that is in turns reflected on a large number of $z$ variables.} In this case the model MT-2IDX, with its smaller memory footprint, is the only viable option to have a good lower bound. It is interesting also to observe that the quality of the heuristic solutions produced by the CP methods is always high, although not always matching the state-of-the-art provided by the nine purely heuristic methods summarized in column HM.

\subsubsection{Instances from Raj et al. \cite{raj2021}} 
A second set of  instances for the MT-PDSTSP was originated from TSPLIB \cite{rei91} in \cite{raj2021}, by defining the missing elements according to what done previously in \cite{mbiadou2018iterative} for the PDSTSP. Only the instances with more than one truck are reported here, since those with one truck boil down the problem to a PDSTSP, and all the optimal solutions can be found in \cite{md23}. The number of customers range between 48 and 229. Full details about the instances can be found in \cite{raj2021}. 
Three methods are considered for comparison purposes: an arc-based MILP model; a branch and price (BP) method based on a set covering MILP model, and a math-heuristic method (MH) based on the same latter model. All the experiments for these methods were run on a computer equipped with an Intel i7-6700 CPU (with 8 cores running at 4 GHz) and 16 GB of RAM, with 3600 seconds for MILP and MH, and 2400 seconds for BC. The CP-based solvers we propose are run for a maximum time of 3600 seconds. The results are summarized in Tables \ref{tr1} and \ref{tr2}.

\begin{table}[h!]													
\caption{MT-PDSVRP - Instances with 2 trucks originated from TSPLIB in Raj et al. \cite{raj2021}.}													
\label{tr1}																	
{	
\hspace{-3cm}				
\footnotesize		
\begin{threeparttable}				
 \begin{tabular}{| l cc| c | c| c | cc|cc|}		\hline	
 \multicolumn{3}{|c|}{Instances}	&		\multicolumn{3}{|c|}{Raj et al.  \cite{raj2021}}	&	\multicolumn{2}{|c|}{MT-3IDX } &	\multicolumn{2}{|c|}{MT-2IDX }	\\
 \multicolumn{3}{|c|}{}						&		\multicolumn{1}{|c|}{ MILP (3600 s)}	&	BP	 (2400 s)&	MH (3600 s)	&	\multicolumn{2}{|c|}{(3600 s)} &	\multicolumn{2}{|c|}{(3600 s)} 	\\
Name		&			$|T|$&			$|D|$		&		UB	&		UB	&	 UB &	LB		&	UB&	LB		&	UB	\\ \hline	
att48\_0\_80	&	2	&	2	&		18206.00				&		17032.00		&	17056.00	&	\textbf{	16940.18	}	&	\textbf{	16940.18}\tnote{*}	&	{	9785.69	}	&	{	17288.00	}	 	\\
att48\_0\_80	&	2	&	4	&		16730.00				&		16500.00		&	16500.00	&	\textbf{	16500.00	}	&		16500.00\tnote{*}	&	{	6562.48	}	&		16564.00		 	\\
att48\_0\_80	&	2	&	6	&		16500.00				&		16500.00		&	16500.00	&	\textbf{	16500.00	}	&		16500.00\tnote{*}	&	{	4980.66	}	&		16500.00		 	\\
berlin52\_0\_80	&	2	&	2	&		3328.20				&		3415.00		&	3415.00	&	\textbf{	3285.00	}	&	\textbf{	3285.00}\tnote{*}	&	{	2240.63	}	&	 	3480.00	 	 	\\
berlin52\_0\_80	&	2	&	4	&		2995.00\tnote{*}		&		2995.00		&	2995.00	&	{	2995.00	}	&		2995.00\tnote{*}	&	{	1498.57	}	&		2995.00		 	\\
berlin52\_0\_80	&	2	&	6	&		2995.00\tnote{*}		&		2995.00		&	2995.00	&	{	2995.00	}	&		2995.00\tnote{*}	&	{	1141.20	}	&		2995.00		 	\\
eil101\_0\_80	&	2	&	2	&		392.00				&		316.00		&	305.40	&	\textbf{	293.00	}	&	\textbf{	293.00}\tnote{*}	&	{	197.44	}	&	{	331.00	}	 	\\
eil101\_0\_80	&	2	&	4	&		408.00				&		266.80		&	253.70	&	\textbf{	232.84	}	&	\textbf{	248.00	}	 	&	{	131.60	}	&	\emph{	283.00	}	 	\\
eil101\_0\_80	&	2	&	6	&		320.00				&		219.00		&	215.30	&	\textbf{	198.00	}	&	\textbf{	212.00	}	 	&	{	99.15	}	&	\emph{	272.00	}	 	\\
gr120\_0\_80	&	2	&	2	&		1145.00				&		806.30		&	764.00	&	\textbf{	703.50	}	&	\textbf{	735.00	}	 	&	{	478.18	}	&	\emph{	866.00	}	 	\\
gr120\_0\_80	&	2	&	4	&		1092.00				&		676.00		&	646.00	&	\textbf{	590.54	}	&	\emph{	676.00	}	 	&	{	318.78	}	&	\emph{	777.00	}	 	\\
gr120\_0\_80	&	2	&	6	&		1082.00				&		769.00		&	581.20	&	\textbf{	518.00	}	&	\textbf{	549.00	}	 	&	{	239.10	}	&	\emph{	688.00	}	 	\\
pr152\_0\_80	&	2	&	2	&		122961.00				&		55477.00		&	48967.00	&	\textbf{	38873.80	}	&	\textbf{	41371.00	}	 	&	{	18792.17	}	&	\emph{	46897.00	}	 	\\
pr152\_0\_80	&	2	&	4	&		-				&		54980.00		&	52353.00	&	\textbf{	33593.10	}	&	\textbf{	41295.58	}	 	&	{	12533.75	}	&	\emph{	44059.00	}	 	\\
pr152\_0\_80	&	2	&	6	&		-				&		50855.40		&	50854.40	&	\textbf{	31020.74	}	&	\textbf{	37413.00	}	 	&	{	9503.63	}	&	\emph{	43573.00	}	 	\\
gr229\_0\_80	&	2	&	2	&		-				&		1403.70		&	1403.70	&	\textbf{	881.91	}	&	\textbf{	1064.97	}	 	&	{	 0.00	}	&	\textbf{	 -	}	 	\\
gr229\_0\_80	&	2	&	4	&		-				&		1346.20		&	1346.20	&	\textbf{	818.33	}	&	\textbf{	1064.97	}	 	&	{	 0.00	}	&	\emph{	 -	}	 	\\
gr229\_0\_80	&	2	&	6	&		-				&		1327.60		&	1327.60	&	\textbf{	776.98	}	&	\textbf{	1056.10	}	 	&	{	 0.00	}	&	\emph{	 -	}	 	\\
\hline
\end{tabular}
\begin{tablenotes}\footnotesize
\item[*] Optimality proven.
\end{tablenotes}
\end{threeparttable}																		
}																	
\end{table}	

The results of Table \ref{tr1} indicates that the CP model MT-3IDX is very suitable on these instances, being able to provide all the state-of-the-art lower bounds and all the best-known heuristic solutions except one, in both cases with many substantial improvements over previously known best results. It is remarkable also that five instances are closed here for the first time. The performance of the model MT-2IDX appear less brilliant, in particular on the larger instances considered, on which no solution or meaningful bound was produced in the given time. 
\begin{table}[h!]												
\caption{MT-PDSVRP - Instances with 3 trucks originated from TSPLIB in Raj et al. \cite{raj2021}.}													
\label{tr2}																	
{			
\hspace{-3cm}						
\footnotesize
\begin{threeparttable}				
 \begin{tabular}{| l cc| c | c| c | cc|cc|}		\hline	
 \multicolumn{3}{|c|}{Instances}	&		\multicolumn{3}{|c|}{Raj et al.  \cite{raj2021}}	&	\multicolumn{2}{|c|}{MT-3IDX } &	\multicolumn{2}{|c|}{MT-2IDX }	\\
 \multicolumn{3}{|c|}{}						&		\multicolumn{1}{|c|}{ MILP (3600 s)}	&	BP	 (2400 s)&	MH (3600 s)	&	\multicolumn{2}{|c|}{(3600 s)} &	\multicolumn{2}{|c|}{(3600 s)} 	\\
Name		&			$|T|$&			$|D|$		&		UB	&		UB	&	 UB &	LB		&	UB&	LB		&	UB	\\ \hline	
att48\_0\_80	&	3	&	2	&		19824.00				&		15218.00		&	14062.00	&	\textbf{	13452.00	}	&	\textbf{	13452.00}\tnote{*}	&	{	7811.15	}	&	\emph{	14652.00	}		\\
att48\_0\_80	&	3	&	4	&		17194.00				&		13394.00		&	14652.00	&	\textbf{	10605.21	}	&	{	13394.00	}		&	{	5690.68	}	&	 	13394.00			\\
att48\_0\_80	&	3	&	6	&		16226.00				&		13756.00		&	13394.00	&	\textbf{	10369.78	}	&		13394.00			&	{	4441.85	}	&		13394.00			\\
berlin52\_0\_80	&	3	&	2	&		3205.00				&		3403.10		&	2995.00	&	\textbf{	2471.35	}	&	\textbf{	2935.00	}		&	{	1777.81	}	&		2995.00			\\
berlin52\_0\_80	&	3	&	4	&		2665.00				&		2635.00		&	2625.00	&	\textbf{	2013.82	}	&		2625.00			&	{	1285.27	}	&	\emph{	2925.00	}		\\
berlin52\_0\_80	&	3	&	6	&		2625.00				&		2625.00		&	2625.00	&	\textbf{	2132.20	}	&		2625.00			&	{	1015.52	}	&	\emph{	2925.00	}		\\
eil101\_0\_80	&	3	&	2	&		-				&		260.00		&	235.00	&	\textbf{	209.00	}	&	\emph{	252.00	}		&	{	157.04	}	&	\emph{	278.00	}		\\
eil101\_0\_80	&	3	&	4	&		-				&		221.00		&	200.10	&	\textbf{	177.00	}	&	\textbf{	198.00	}		&	{	112.23	}	&	\emph{	227.00	}		\\
eil101\_0\_80	&	3	&	6	&		-				&		228.00		&	181.00	&	\textbf{	155.20	}	&	\emph{	218.00	}		&	{	87.91	}	&	\emph{	204.00	}		\\
gr120\_0\_80	&	3	&	2	&		-				&		611.00		&	597.00	&	\textbf{	509.50	}	&	\emph{	688.00	}		&	{	381.29	}	&	\emph{	734.00	}		\\
gr120\_0\_80	&	3	&	4	&		-				&		536.00		&	528.00	&	\textbf{	444.00	}	&	\emph{	688.00	}		&	{	272.58	}	&	\emph{	601.00	}		\\
gr120\_0\_80	&	3	&	6	&		-				&		590.90		&	527.20	&	\textbf{	397.50	}	&	\textbf{	526.00	}		&	{	212.84	}	&	\emph{	547.00	}		\\
pr152\_0\_80	&	3	&	2	&		-				&		48135.00		&	48135.00	&	\textbf{	25154.41	}	&	{	38395.00	}		&	{	15033.25	}	&	\textbf{	36443.00	}		\\
pr152\_0\_80	&	3	&	4	&		-				&		43024.00		&	43024.00	&	\textbf{	21913.95	}	&	\textbf{	36940.00	}		&	{	10799.23	}	&	{	37183.00	}		\\
pr152\_0\_80	&	3	&	6	&		-				&		46403.00		&	46403.00	&	\textbf{	20584.69	}	&	{	36940.00	}		&	{	8526.30	}	&	\textbf{	36281.00	}		\\
gr229\_0\_80	&	3	&	2	&		-				&		1145.40		&	1145.40	&	\textbf{	579.62	}	&	\textbf{	981.59	}		&	{	 0.00	}	&	\emph{	 -	}		\\
gr229\_0\_80	&	3	&	4	&		-				&		1132.30		&	1132.30	&	\textbf{	561.47	}	&	\textbf{	969.92	}		&	{	 0.00	}	&	\emph{	 -	}		\\
gr229\_0\_80	&	3	&	6	&		-				&		1130.00		&	1130.00	&	\textbf{	517.37	}	&	\textbf{	1004.89	}		&	{	 0.00	}	&	\emph{	 -	}		\\
\hline
\end{tabular}
\begin{tablenotes}\footnotesize
\item[*] Optimality proven.
\end{tablenotes}
\end{threeparttable}																		
}																	
\end{table}

The results of Table \ref{tr2} confirm the previous impressions. Only, in this case the CP model MT-3IDX otperforms the other methods less strongly. This indicates that its performance degrade when the number of trucks considered increases: On a few instances it found sub-optimal heuristic solutions, although it remarkably provided an optimality proof for the first instance. Note that the model MT-2IDX is now able to produce two new state-of-the-art heuristic solutions, indicating that this latter method might be indicated for problems with more trucks available.

\subsection{MC-PDSVRP} \label{rvrpc}	

The instances considered in this section are those originally proposed in Nguyen et al. \cite{nguyen2022}. They are based on several statistical data shared by logistics service providers and common practices, so although being generated as random, they can be considered as realistic. The numbers of customers considered are 30, 50, 100, 200 and 400 (the first element of the name of each instance reflects this information). We refer the reader to \cite{nguyen2022} for the full details of the instances ranging from the number of trucks and drones, to loading capacities, speeds, battery endurance for the drones, etc. 

The results of our experiments are summarized in Tables \ref{tv30}, \ref{tv50}, \ref{tv100}, \ref{tv200} and \ref{tv400}, divided according the size of the instances. On top of the CP models discussed in Section \ref{cpc}, that are run for a maximum of 1800 seconds, the methods considered in the comparison are those presented in \cite{nguyen2022}, namely a MILP model and \emph{Ruin and Recreate} heuristic algorithm (RR). The experiments for the latter methods were run on a computer equipped with an AMD Ryzen 3700X CPU running at 4.0 GHz and 16 GB of RAM, with a time limit of 3600 seconds. Note that the RR method was run 30 times, and the best result is reported here.

\begin{table}[h!]													
\caption{MC-PDSVRP - Instances with 30 customers from Nguyen et al. \cite{nguyen2022}.}													
\label{tv30}																	
{
\hspace{-.5cm}							
\footnotesize		
\begin{threeparttable}		
 \begin{tabular}{| l | c | c| cc | cc|}		\hline	
Instances	&		\multicolumn{2}{|c|}{Nguyen et al. \cite{nguyen2022}}	&		\multicolumn{2}{|c|}{MC-3IDX}&\multicolumn{2}{|c|}{MC-2IDX}\\
	&		         \multicolumn{1}{|c|}{MILP	(3600 sec)}	&	\multicolumn{1}{|c|}{RR (30x3600 sec)}&	\multicolumn{2}{|c}{(1800 sec)}	 &	\multicolumn{2}{|c|}{(1800 sec)}	 		\\
				&		UB		&				UB	&				LB		&UB	 &		LB		&UB	 	\\ \hline
30-r-0-c	&	128.67\tnote{*}	&	128.67	&	{	101.69	}	&	{	128.67	}		&	\emph{	103.11	}	&	{	128.67	}		\\
30-r-0-e	&	300.85		&	300.85	&	{	104.85	}	&	{	300.85	}		&	\textbf{	150.61	}	&	{	300.85	}		\\
30-r-0-r	&	199.76		&	199.76	&	{	110.71	}	&	{	199.76	}		&	\textbf{	125.00	}	&	{	199.76	}		\\
30-r-1-c	&	58.30\tnote{*}	&	58.30	&	{	58.30	}	&	{	58.30}\tnote{*}	&	{	58.30	}	&	{	58.30}\tnote{*}	\\
30-r-1-e	&	243.11		&	243.11	&	{	117.48	}	&	{	243.11	}		&	\textbf{	148.43	}	&	{	243.11	}		\\
30-r-1-r	&	164.75		&	164.75	&	{	106.82	}	&	{	164.75	}		&	\textbf{	116.50	}	&	{	164.75	}		\\
30-r-2-c	&	91.57\tnote{*}	&	91.57	&	{	69.04	}	&	{	91.57	}		&	\emph{	80.43	}	&	{	91.57	}		\\
30-r-2-e	&	235.25		&	235.25	&	{	111.62	}	&	{	235.25	}		&	\textbf{	151.10	}	&	{	235.25	}		\\
30-r-2-r	&	211.26		&	211.26	&	{	101.11	}	&	{	211.26	}		&	\textbf{	135.88	}	&	{	211.26	}		\\
30-c-0-c	&	82.15\tnote{*}	&	82.15	&	\emph{	60.08	}	&	{	82.15	}		&	{	82.15	}	&	{	82.15}\tnote{*}	\\
30-c-0-e	&	150.44\tnote{*}	&	150.44	&	\emph{	103.78	}	&	{	150.44	}		&	\emph{	110.23	}	&	{	150.44	}		\\
30-c-0-r	&	112.72\tnote{*}	&	112.72	&	\emph{	83.08	}	&	{	112.72	}		&	\emph{	91.44	}	&	{	112.72	}		\\
30-c-1-c	&	39.98\tnote{*}	&	39.98	&	{	39.98	}	&	{	39.98}\tnote{*}	&	{	39.98	}	&	{	39.98}\tnote{*}	\\
30-c-1-e	&	172.49		&	172.49	&	\textbf{	97.03	}	&	{	172.49	}		&	{	83.94	}	&	{	172.49	}		\\
30-c-1-r	&	61.28\tnote{*}	&	61.28	&	{	61.28	}	&	{	61.28}\tnote{*}	&	{	61.28	}	&	{	61.28}\tnote{*}	\\
30-c-2-c	&	77.00\tnote{*}	&	77.00	&	{	77.00	}	&	{	77.00}\tnote{*}	&	{	77.00	}	&	{	77.00}\tnote{*}	\\
30-c-2-e	&	166.82		&	166.82	&	\textbf{	98.49	}	&	{	166.82	}		&	{	95.89	}	&	{	166.82	}		\\
30-c-2-r	&	88.81\tnote{*}	&	88.81	&	\emph{	59.65	}	&	{	88.81	}		&	\emph{	80.22	}	&	{	88.81	}		\\
30-rc-0-c	&	81.06\tnote{*}	&	81.06	&	{	81.06	}	&	{	81.06}\tnote{*}	&	\emph{	74.47	}	&	{	81.06	}		\\
30-rc-0-e	&	181.89		&	181.89	&	{	111.08	}	&	{	181.89	}		&	\textbf{	115.78	}	&	{	181.89	}		\\
30-rc-0-r	&	146.73		&	146.73	&	{	88.40	}	&	{	146.73	}		&	\textbf{	98.85	}	&	{	146.73	}		\\
30-rc-1-c	&	56.48\tnote{*}	&	56.48	&	{	56.48	}	&	{	56.48}\tnote{*}	&	{	56.48	}	&	{	56.48}\tnote{*}	\\
30-rc-1-e	&	176.57		&	176.57	&	{	115.70	}	&	{	176.57	}		&	\textbf{	118.09	}	&	{	176.57	}		\\
30-rc-1-r	&	97.91\tnote{*}	&	97.91	&	{	97.91	}	&	{	97.91}\tnote{*}	&	\emph{	71.96	}	&	{	97.91	}		\\
30-rc-2-c	&	77.00\tnote{*}	&	77.00	&	\emph{	63.74	}	&	{	77.00	}		&	{	77.00	}	&	{	77.00}\tnote{*}	\\
30-rc-2-e	&	187.18		&	187.18	&	{	71.50	}	&	{	187.18	}		&	\textbf{	104.96	}	&	{	187.18	}		\\
30-rc-2-r	&	94.21\tnote{*}	&	94.21	&	\emph{	77.67	}	&	{	94.21	}		&	\emph{	86.67	}	&	{	94.21	}		\\
\hline
\end{tabular}	
\begin{tablenotes}\footnotesize
\item[*] Optimality proven.
\end{tablenotes}
\end{threeparttable}																	
}																	
\end{table}

\begin{table}														
\caption{MC-PDSVRP - Instances with 50 customers from Nguyen et al. \cite{nguyen2022}.}													
\label{tv50}																	
{	
\hspace{-.5cm}						
\footnotesize		
\begin{threeparttable}		
 \begin{tabular}{| l | c | c| cc | cc|}		\hline	
Instances	&		\multicolumn{2}{|c|}{Nguyen et al. \cite{nguyen2022}}	&		\multicolumn{2}{|c|}{MC-3IDX}&\multicolumn{2}{|c|}{MC-2IDX}\\
	&		         \multicolumn{1}{|c|}{MILP	(3600 sec)}	&	\multicolumn{1}{|c|}{RR (30x3600 sec)}&	\multicolumn{2}{|c}{(1800 sec)}	 &	\multicolumn{2}{|c|}{(1800 sec)}	 		\\
				&		UB		&				UB	&				LB		&UB	 &		LB		&UB	 	\\ \hline
50-r-0-c	&	161.13\tnote{*}	&	161.13	&	\emph{	110.62	}	&	{	161.13	}		&	\emph{	131.27	}	&	\emph{	161.56	}		\\
50-r-0-e	&	-		&	310.50	&	{	129.68	}	&	{	310.50	}		&	\textbf{	151.57	}	&	{	310.50	}		\\
50-r-0-r	&	-		&	246.53	&	{	110.02	}	&	{	246.53	}		&	\textbf{	137.85	}	&	{	246.53	}		\\
50-r-1-c	&	177.21		&	177.21	&	{	103.69	}	&	{	177.21	}		&	\textbf{	135.27	}	&	{	177.21	}		\\
50-r-1-e	&	216.97		&	216.97	&	{	129.97	}	&	{	216.97	}		&	\textbf{	153.12	}	&	{	216.97	}		\\
50-r-1-r	&	-		&	228.11	&	{	116.45	}	&	{	228.11	}		&	\textbf{	148.57	}	&	{	228.11	}		\\
50-r-2-c	&	182.09		&	181.49	&	{	108.59	}	&	{	182.62	}		&	\textbf{	145.83	}	&	\emph{	182.62	}		\\
50-r-2-e	&	-		&	375.05	&	{	127.88	}	&	{	375.05	}		&	\textbf{	153.27	}	&	{	375.05	}		\\
50-r-2-r	&	-		&	313.32	&	{	109.95	}	&	{	313.32	}		&	\textbf{	162.11	}	&	{	313.32	}		\\
50-c-0-c	&	89.04		&	89.04	&	{	47.73	}	&	{	89.04	}		&	\textbf{	78.64	}	&	{	89.04	}		\\
50-c-0-e	&	183.26		&	182.91	&	{	78.01	}	&	{	182.91	}		&	\textbf{	105.82	}	&	{	182.91	}		\\
50-c-0-r	&	112.17		&	112.17	&	{	85.09	}	&	{	112.17	}		&	\textbf{	89.09	}	&	{	112.17	}		\\
50-c-1-c	&	81.76\tnote{*}	&	81.76	&	\emph{	39.38	}	&	{	81.76	}		&	\emph{	75.97	}	&	{	81.76	}		\\
50-c-1-e	&	133.57		&	133.57	&	{	81.44	}	&	{	133.57	}		&	\textbf{	92.35	}	&	{	133.57	}		\\
50-c-1-r	&	149.24		&	148.86	&	{	56.78	}	&	\emph{	149.06	}		&	\textbf{	83.95	}	&	\emph{	149.06	}		\\
50-c-2-c	&	77.16\tnote{*}	&	77.16	&	\emph{	54.26	}	&	{	77.16	}		&	\emph{	75.31	}	&	{	77.16	}		\\
50-c-2-e	&	214.25		&	214.25	&	{	74.24	}	&	\emph{	214.66	}		&	\textbf{	112.86	}	&	\emph{	214.66	}		\\
50-c-2-r	&	136.60		&	136.51	&	{	67.30	}	&	{	136.51	}		&	\textbf{	90.66	}	&	{	136.51	}		\\
50-rc-0-c	&	116.80\tnote{*}	&	116.80	&	\emph{	70.01	}	&	{	116.80	}		&	\emph{	107.92	}	&	{	116.80	}		\\
50-rc-0-e	&	-		&	242.02	&	{	74.24	}	&	{	242.02	}		&	\textbf{	124.28	}	&	{	242.02	}		\\
50-rc-0-r	&	217.87		&	217.87	&	{	95.88	}	&	{	217.87	}		&	\textbf{	109.28	}	&	{	217.87	}		\\
50-rc-1-c	&	110.44\tnote{*}	&	110.44	&	\emph{	82.98	}	&	{	110.44	}		&	\emph{	104.16	}	&	{	110.44	}		\\
50-rc-1-e	&	201.62		&	201.62	&	{	112.97	}	&	{	201.62	}		&	\textbf{	124.27	}	&	{	201.62	}		\\
50-rc-1-r	&	173.16		&	173.16	&	{	90.57	}	&	\emph{	173.61	}		&	\textbf{	116.16	}	&	\emph{	173.61	}		\\
50-rc-2-c	&	129.82		&	129.65	&	{	64.33	}	&	{	129.82	}		&	\textbf{	104.74	}	&	{	129.65	}		\\
50-rc-2-e	&	-		&	383.49	&	{	109.51	}	&	{	383.49	}		&	\textbf{	129.50	}	&	{	383.49	}		\\
50-rc-2-r	&	197.20		&	186.31	&	{	108.11	}	&	{	196.97	}		&	\textbf{	108.82	}	&	{	186.31	}		\\
\hline
\end{tabular}		
\begin{tablenotes}\footnotesize
\item[*] Optimality proven.
\end{tablenotes}
\end{threeparttable}																
}																	
\end{table}	

The results reported in Tables \ref{tv30} and \ref{tv50} suggest that the MILP model is the most viable method to solve these small problems with 30 and 50 customers: several instances are solved to optimality within the given time. However, on some of the instances with 50 customers, MILP is not able to provide any solution. The RR heuristic also obtains remarkably good results, always matching or improving the best known figures. We believe that taking the best of 30 runs was an important factor in these results, because this allows a good random exploration of the search space. The CP models perform reasonably well on these instances in terms of heuristic solution found, and are able to provide good lower bounds (no lower bound was provided in previous publications). It is worth noticing that the model MC-2IDX seems to perform better than MC-3 IDX on these \emph{Minimum Cost} instances, especially on the instances with 50 customers. This revert what was observed in Section \ref{rvrpt} for the \emph{Minimum Time} instances, and can be explained by the completely different objective function between the two problems considered, and their general characteristics.

\begin{table}[h!]													
\caption{MC-PDSVRP - Instances with 100 customers from Nguyen et al. \cite{nguyen2022}.}													
\label{tv100}
{	
\hspace{-.5cm}						
\footnotesize				
 \begin{tabular}{| l | c | c| cc | cc|}		\hline	
Instances	&		\multicolumn{2}{|c|}{Nguyen et al. \cite{nguyen2022}}	&		\multicolumn{2}{|c|}{MC-3IDX}&\multicolumn{2}{|c|}{MC-2IDX}\\
	&		         \multicolumn{1}{|c|}{MILP	(3600 sec)}	&	\multicolumn{1}{|c|}{RR (30x3600 sec)}&	\multicolumn{2}{|c}{(1800 sec)}	 &	\multicolumn{2}{|c|}{(1800 sec)}	 		\\
				&		UB		&				UB	&				LB		&UB	 &		LB		&UB	 	\\ \hline
100-r-0-c	&	 -		&	272.77	&	{	143.76	}	&	\emph{	 -	}		&	\textbf{	222.27	}	&	\emph{	277.90	}		\\
100-r-0-e	&	 -		&	506.07	&	{	206.21	}	&	\emph{	 -	}		&	\textbf{	250.30	}	&	\emph{	 -	}		\\
100-r-0-r	&	 -		&	325.29	&	{	140.62	}	&	\emph{	 -	}		&	\textbf{	226.44	}	&	\emph{	 -	}		\\
100-r-1-c	&	 -		&	302.46	&	{	167.13	}	&	\emph{	 -	}		&	\textbf{	240.03	}	&	\emph{	 -	}		\\
100-r-1-e	&	 -		&	470.06	&	{	198.41	}	&	\emph{	 -	}		&	\textbf{	251.18	}	&	\emph{	 -	}		\\
100-r-1-r	&	 -		&	350.97	&	{	179.63	}	&	\emph{	 -	}		&	\textbf{	230.43	}	&	\emph{	 -	}		\\
100-r-2-c	&	 -		&	284.77	&	{	169.39	}	&	\emph{	288.00	}		&	\textbf{	232.12	}	&	\emph{	 -	}		\\
100-r-2-e	&	 -		&	435.91	&	{	208.16	}	&	\emph{	 -	}		&	\textbf{	265.24	}	&	\emph{	 -	}		\\
100-r-2-r	&	 -		&	320.98	&	{	162.92	}	&	\emph{	361.47	}		&	\textbf{	226.69	}	&	\emph{	 -	}		\\
100-c-0-c	&	 -		&	218.59	&	{	92.16	}	&	\emph{	221.59	}		&	\textbf{	170.73	}	&	\emph{	 -	}		\\
100-c-0-e	&	 -		&	315.84	&	{	93.02	}	&	\emph{	 -	}		&	\textbf{	175.48	}	&	\emph{	 -	}		\\
100-c-0-r	&	 -		&	252.39	&	{	81.02	}	&	\emph{	276.95	}		&	\textbf{	161.81	}	&	\emph{	 -	}		\\
100-c-1-c	&	 -		&	154.61	&	{	87.41	}	&	\emph{	179.75	}		&	\textbf{	112.94	}	&	\emph{	 -	}		\\
100-c-1-e	&	 -		&	222.15	&	{	103.33	}	&	\emph{	236.44	}		&	\textbf{	119.22	}	&	\emph{	 -	}		\\
100-c-1-r	&	 -		&	110.95	&	{	64.71	}	&	{	110.95	}		&	\textbf{	98.77	}	&	{	110.95	}		\\
100-c-2-c	&	 -		&	108.35	&	{	73.34	}	&	{	108.35	}		&	\textbf{	95.81	}	&	\emph{	 -	}		\\
100-c-2-e	&	 -		&	118.71	&	{	70.71	}	&	\emph{	132.28	}		&	\textbf{	103.69	}	&	\emph{	 -	}		\\
100-c-2-r	&	 -		&	86.28	&	{	49.08	}	&	\emph{	102.47	}		&	\textbf{	86.00	}	&	\emph{	88.01	}		\\
100-rc-0-c	&	 -		&	252.37	&	{	128.09	}	&	\emph{	254.03	}		&	\textbf{	199.69	}	&	\emph{	 -	}		\\
100-rc-0-e	&	 -		&	472.95	&	{	168.08	}	&	\emph{	 -	}		&	\textbf{	233.41	}	&	\emph{	 -	}		\\
100-rc-0-r	&	 -		&	337.61	&	{	157.65	}	&	\emph{	339.09	}		&	\textbf{	205.28	}	&	\emph{	 -	}		\\
100-rc-1-c	&	 -		&	276.78	&	{	152.69	}	&	\emph{	305.51	}		&	\textbf{	209.72	}	&	\emph{	 -	}		\\
100-rc-1-e	&	 -		&	432.44	&	{	190.81	}	&	\emph{	490.66	}		&	\textbf{	246.80	}	&	\emph{	 -	}		\\
100-rc-1-r	&	 -		&	423.33	&	{	163.48	}	&	\emph{	426.76	}		&	\textbf{	220.13	}	&	\emph{	 -	}		\\
100-rc-2-c	&	 -		&	260.45	&	{	125.31	}	&	\emph{	286.73	}		&	\textbf{	198.66	}	&	\emph{	 -	}		\\
100-rc-2-e	&	 -		&	358.24	&	{	174.37	}	&	\emph{	406.12	}		&	\textbf{	228.51	}	&	\emph{	 -	}		\\
100-rc-2-r	&	 -		&	331.97	&	{	138.87	}	&	\emph{	345.25	}		&	\textbf{	207.27	}	&	\emph{	 -	}		\\
\hline
\end{tabular}																	
}																	
\end{table}

\begin{table}[h!]														
\caption{MC-PDSVRP - Instances with 200 customers from Nguyen et al. \cite{nguyen2022}.}													
\label{tv200}																	
{	
\hspace{-0cm}						
\footnotesize				
 \begin{tabular}{| l | c | c| cc | cc|}		\hline	
Instances	&		\multicolumn{2}{|c|}{Nguyen et al. \cite{nguyen2022}}	&		\multicolumn{2}{|c|}{MC-3IDX}&\multicolumn{2}{|c|}{MC-2IDX}\\
	&		         \multicolumn{1}{|c|}{MILP	(3600 sec)}	&	\multicolumn{1}{|c|}{RR (30x3600 sec)}&	\multicolumn{2}{|c}{(1800 sec)}	 &	\multicolumn{2}{|c|}{(1800 sec)}	 		\\
				&		UB		&				UB	&				LB		&UB	 &		LB		&UB	 	\\ \hline
200-r-0-c	&	 -		&	430.42	&	{	134.99	}	&	\emph{	 -	}		&	\textbf{	353.51	}	&	{	 -	}		\\
200-r-0-e	&	 -		&	728.00	&	{	0.00	}	&	\emph{	 -	}		&	\textbf{	368.88	}	&	{	 -	}		\\
200-r-0-r	&	 -		&	528.29	&	{	102.71	}	&	\emph{	 -	}		&	\textbf{	355.73	}	&	{	 -	}		\\
200-r-1-c	&	 -		&	434.66	&	{	228.69	}	&	\emph{	 -	}		&	\textbf{	348.04	}	&	{	 -	}		\\
200-r-1-e	&	 -		&	601.06	&	{	0.13	}	&	\emph{	 -	}		&	\textbf{	369.88	}	&	{	 -	}		\\
200-r-1-r	&	 -		&	481.83	&	{	128.30	}	&	\emph{	 -	}		&	\textbf{	349.86	}	&	{	 -	}		\\
200-r-2-c	&	 -		&	446.15	&	{	0.00	}	&	\emph{	 -	}		&	\textbf{	357.04	}	&	{	 -	}		\\
200-r-2-e	&	 -		&	606.05	&	{	0.00	}	&	\emph{	 -	}		&	\textbf{	374.19	}	&	{	 -	}		\\
200-r-2-r	&	 -		&	497.72	&	{	0.00	}	&	\emph{	 -	}		&	\textbf{	350.49	}	&	{	 -	}		\\
200-c-0-c	&	 -		&	234.40	&	{	0.00	}	&	\emph{	 -	}		&	\textbf{	184.05	}	&	{	 -	}		\\
200-c-0-e	&	 -		&	442.74	&	{	0.00	}	&	\emph{	 -	}		&	\textbf{	189.35	}	&	{	 -	}		\\
200-c-0-r	&	 -		&	270.76	&	{	0.00	}	&	\emph{	 -	}		&	\textbf{	173.94	}	&	{	 -	}		\\
200-c-1-c	&	 -		&	225.55	&	{	78.21	}	&	\emph{	 -	}		&	\textbf{	186.77	}	&	{	 -	}		\\
200-c-1-e	&	 -		&	461.57	&	{	0.00	}	&	\emph{	 -	}		&	\textbf{	211.23	}	&	{	 -	}		\\
200-c-1-r	&	 -		&	239.91	&	{	0.00	}	&	\emph{	 -	}		&	\textbf{	198.22	}	&	{	 -	}		\\
200-c-2-c	&	 -		&	263.26	&	{	0.00	}	&	\emph{	 -	}		&	\textbf{	193.47	}	&	{	 -	}		\\
200-c-2-e	&	 -		&	486.05	&	{	0.00	}	&	\emph{	 -	}		&	\textbf{	212.55	}	&	{	 -	}		\\
200-c-2-r	&	 -		&	322.61	&	{	0.00	}	&	\emph{	 -	}		&	\textbf{	189.90	}	&	{	 -	}		\\
200-rc-0-c	&	 -		&	354.23	&	{	0.00	}	&	\emph{	 -	}		&	\textbf{	280.87	}	&	{	 -	}		\\
200-rc-0-e	&	 -		&	575.94	&	{	0.00	}	&	\emph{	 -	}		&	\textbf{	297.83	}	&	{	 -	}		\\
200-rc-0-r	&	 -		&	484.66	&	{	0.00	}	&	\emph{	 -	}		&	\textbf{	293.68	}	&	{	 -	}		\\
200-rc-1-c	&	 -		&	370.01	&	{	0.00	}	&	\emph{	 -	}		&	\textbf{	280.68	}	&	{	 -	}		\\
200-rc-1-e	&	 -		&	591.58	&	{	0.00	}	&	\emph{	 -	}		&	\textbf{	296.25	}	&	{	 -	}		\\
200-rc-1-r	&	 -		&	411.34	&	{	0.00	}	&	\emph{	 -	}		&	\textbf{	267.01	}	&	{	 -	}		\\
200-rc-2-c	&	 -		&	372.30	&	{	0.00	}	&	\emph{	 -	}		&	\textbf{	299.35	}	&	{	 -	}		\\
200-rc-2-e	&	 -		&	534.45	&	{	0.00	}	&	\emph{	 -	}		&	\textbf{	330.13	}	&	{	 -	}		\\
200-rc-2-r	&	 -		&	432.17	&	{	0.00	}	&	\emph{	 -	}		&	\textbf{	299.16	}	&	{	 -	}		\\
\hline
\end{tabular}																	
}																	
\end{table}	

\begin{table}[h!]													
\caption{MC-PDSVRP - Instances with 400 customers from Nguyen et al. \cite{nguyen2022}.}													
\label{tv400}
{			
\hspace{-0cm}		
\footnotesize				
 \begin{tabular}{| l | c | c| cc | cc|}		\hline	
Instances	&		\multicolumn{2}{|c|}{Nguyen et al. \cite{nguyen2022}}	&		\multicolumn{2}{|c|}{MC-3IDX}&\multicolumn{2}{|c|}{MC-2IDX}\\
	&		         \multicolumn{1}{|c|}{MILP	(3600 sec)}	&	\multicolumn{1}{|c|}{RR (30x3600 sec)}&	\multicolumn{2}{|c}{(1800 sec)}	 &	\multicolumn{2}{|c|}{(1800 sec)}	 		\\
				&		UB		&				UB	&				LB		&UB	 &		LB		&UB	 	\\ \hline
400-r-0-c	&	 -		&	604.94	&	{	0.00	}	&	\emph{	 -	}		&	\textbf{	500.91	}	&	\emph{	 -	}		\\
400-r-0-e	&	 -		&	801.45	&	{	0.00	}	&	\emph{	 -	}		&	\textbf{	518.74	}	&	\emph{	 -	}		\\
400-r-0-r	&	 -		&	626.46	&	{	0.00	}	&	\emph{	 -	}		&	\textbf{	518.62	}	&	\emph{	 -	}		\\
400-r-1-c	&	 -		&	578.77	&	{	0.00	}	&	\emph{	 -	}		&	\textbf{	487.76	}	&	\emph{	 -	}		\\
400-r-1-e	&	 -		&	772.30	&	{	0.00	}	&	\emph{	 -	}		&	\textbf{	445.80	}	&	\emph{	 -	}		\\
400-r-1-r	&	 -		&	584.07	&	{	0.00	}	&	\emph{	 -	}		&	\textbf{	477.13	}	&	\emph{	 -	}		\\
400-r-2-c	&	 -		&	587.88	&	{	0.00	}	&	\emph{	 -	}		&	\textbf{	483.98	}	&	\emph{	 -	}		\\
400-r-2-e	&	 -		&	771.16	&	{	0.00	}	&	\emph{	 -	}		&	\textbf{	510.02	}	&	\emph{	 -	}		\\
400-r-2-r	&	 -		&	692.91	&	{	0.00	}	&	\emph{	 -	}		&	\textbf{	454.20	}	&	\emph{	 -	}		\\
400-c-0-c	&	 -		&	271.82	&	{	0.00	}	&	\emph{	 -	}		&	\textbf{	207.27	}	&	\emph{	 -	}		\\
400-c-0-e	&	 -		&	378.98	&	{	0.00	}	&	\emph{	 -	}		&	\textbf{	253.55	}	&	\emph{	 -	}		\\
400-c-0-r	&	 -		&	403.85	&	{	0.00	}	&	\emph{	 -	}		&	\textbf{	248.69	}	&	\emph{	 -	}		\\
400-c-1-c	&	 -		&	187.20	&	{	0.00	}	&	\emph{	 -	}		&	\textbf{	136.34	}	&	\emph{	 -	}		\\
400-c-1-e	&	 -		&	249.01	&	{	0.00	}	&	\emph{	 -	}		&	\textbf{	142.46	}	&	\emph{	 -	}		\\
400-c-1-r	&	 -		&	219.23	&	{	0.00	}	&	\emph{	 -	}		&	\textbf{	146.12	}	&	\emph{	 -	}		\\
400-c-2-c	&	 -		&	152.79	&	{	0.00	}	&	\emph{	 -	}		&	\textbf{	141.57	}	&	\emph{	 -	}		\\
400-c-2-e	&	 -		&	245.54	&	{	0.00	}	&	\emph{	 -	}		&	\textbf{	145.22	}	&	\emph{	 -	}		\\
400-c-2-r	&	 -		&	231.89	&	{	0.00	}	&	\emph{	 -	}		&	\textbf{	154.86	}	&	\emph{	 -	}		\\
400-rc-0-c	&	 -		&	497.90	&	{	0.00	}	&	\emph{	 -	}		&	\textbf{	420.36	}	&	\emph{	 -	}		\\
400-rc-0-e	&	 -		&	700.96	&	{	0.00	}	&	\emph{	 -	}		&	\textbf{	456.69	}	&	\emph{	 -	}		\\
400-rc-0-r	&	 -		&	604.38	&	{	0.00	}	&	\emph{	 -	}		&	\textbf{	437.29	}	&	\emph{	 -	}		\\
400-rc-1-c	&	 -		&	518.40	&	{	0.00	}	&	\emph{	 -	}		&	\textbf{	395.26	}	&	\emph{	 -	}		\\
400-rc-1-e	&	 -		&	784.02	&	{	0.00	}	&	\emph{	 -	}		&	\textbf{	424.46	}	&	\emph{	 -	}		\\
400-rc-1-r	&	 -		&	531.44	&	{	0.00	}	&	\emph{	 -	}		&	\textbf{	396.43	}	&	\emph{	 -	}		\\
400-rc-2-c	&	 -		&	510.64	&	{	0.00	}	&	\emph{	 -	}		&	\textbf{	424.09	}	&	\emph{	 -	}		\\
400-rc-2-e	&	 -		&	723.63	&	{	0.00	}	&	\emph{	 -	}		&	\textbf{	449.39	}	&	\emph{	 -	}		\\
400-rc-2-r	&	 -		&	625.31	&	{	0.00	}	&	\emph{	 -	}		&	\textbf{	448.43	}	&	\emph{	 -	}		\\
\hline
\end{tabular}																	
}																	
\end{table}	

The results summarized in Tables \ref{tv100}, \ref{tv200} and \ref{tv400} do not include the results of the MILP model, since these larger instances with 100, 200 and 400 customers were out of reach for it. The RR heuristic produce good quality solutions, as we are able to certify here with the first lower bounds ever produced on these instances. On the other hand, the relatively small gaps between upper and lower bounds on some of the instances, also suggest that the quality provided by the CP models are good. In details, MC-2IDX scales up very well while MC-3IDX is not able to produce significant bounds, due to the very large number of $z$ variables necessary to model the problem. Unfortunately these good performance of the model MC-2IDX in terms of lower bounds does not get reflected on the quality of the heuristic solutions: very often the method is not capable of producing any feasible solution at all in the given computation time.
										
\section{Conclusions} \label{conc}	
In this paper we have discussed some different Constraint Programming models to describe two versions of the Parallel Drone Scheduling Vehicle Routing Problem that were recently proposed in the literature. 

Experimental results suggest that solving these models can lead to many improved state-of-the-art results. In particular, the new models seem to provide the new reference for producing high quality lower bounds on the optimal solution costs, but they are also able to produce several new best-known heuristic solutions, and even to close for the first time several of the instances considered. 

In our opinion the flourishing literature on specialization of Parallel Drone Scheduling Vehicle Routing Problems, aiming at introducing more and more realistic aspects could be greatly benefit from our contributions, given the high flexibility provided by Constraint Programming, and the assumption that the results we have obtained here could be replicated. This is material for future work.
			
\section*{Acknowledgements}
The authors are grateful to Prof. Ho\`{a}ng Ha Minh for the useful discussions and suggestions.

\section*{Declaration of competing interest}
The authors declare that they have no known competing financial interests or personal
relationships that could have appeared to influence the work reported in this paper.
\bibliographystyle{plain}

\end{document}